\documentclass{ifacconf}

\usepackage{graphicx}      
\usepackage{natbib}        
\usepackage{color} 
\usepackage{tabularx}
\usepackage{latexsym,amssymb,amsmath,amsfonts,amscd} 
\usepackage[all]{xy}
\usepackage{mathrsfs} 
\usepackage{tikz} 
\usetikzlibrary{shapes,patterns,arrows,decorations.pathreplacing}
\usepackage{epigraph,enumerate}
\usepackage{stmaryrd}
\usepackage{tikz-cd}
\usepackage{mathtools}

\usepackage{latexsym, amssymb, amsmath}
\usepackage{enumerate}
\usepackage{flushend} 
\usepackage{mathrsfs} 
\usepackage{stfloats}
\def\custombibliography#1{
 \normalsize
\section*{\centering References}
 \list
 {[\arabic{enumi}]}{\settowidth\labelwidth{[#1]}\leftmargin\labelwidth
 \setlength{\itemsep}{.1em}
 \advance\leftmargin\labelsep
 \usecounter{enumi}}
 \def\newblock{\hskip .11em plus .33em minus -.07em}
 \sloppy
 \sfcode`\.=1000\relax}

\def\L2{{\cal L}_2}

\newcommand\bull{\vrule height .9ex width .8ex depth -.1ex } 

\newcommand\re{\rm I\! R}
\newcommand\cdcout[1]{} 

\newcommand{\rv}[1]{\boldsymbol{#1}} 
\newcommand{\RomanNumber}[1]{\uppercase\expandafter{\romannumeral #1}}
\newcommand{\romannumber}[1]{\lowercase\expandafter{\romannumeral #1}}

\DeclareMathAlphabet{\mathpzc}{OT1}{pzc}{m}{it}

\def\1{\rv 1} 

\usepackage{flushend} 
\usepackage{stfloats}
\usepackage{color} 
\usepackage{latexsym,amssymb,amsmath} 
\usepackage{mathrsfs,mathtools} 
\usepackage{multicol,enumerate}
\usepackage{subcaption}

\def\allpoly{\mbox{$\re\langle X \rangle$}}

\def\allpolyx0degn{\mbox{$P_n$}}

\def\allwords{\mbox{$X^{\ast}$}}

\def\allseries{\mbox{$\re\langle\langle X \rangle\rangle$}}
\def\allseries#1{\mbox{$\re^{#1}\langle\langle X \rangle\rangle$}}

\def\allproperseries#1{\mbox{$\re_{p}^{#1}\, \langle\langle X \rangle\rangle$}}

\def\allseriesX1{\mbox{$\re [[ X_1 ]]$}}

\def\bull{\rule{0.08in}{0.08in}} 

\newcommand{\comment}[1]{} 

\def\eqref#1{(\ref{#1})} 

\def\mbf#1{\hbox{\mathversion{bold}$#1$}} 
\def\modcomp{\:\tilde{\circ}\,} 
\def\norm#1{\Vert#1\Vert}

\def\openbull{\framebox[0.08in][c]{$\;$}} 

\def\re{{\mathbb R}} 

\def\shuffle{{\scriptscriptstyle \;\sqcup \hspace*{-0.05cm}\sqcup\;}}

\def\vec{{\rm vec}}

\def\begals{\[\begin{aligned}}
\def\endals{\end{aligned}\]}
\def\begal{\begin{align*}}
\def\endal{\end{align*}}
\def\begce{\begin{center}}
\def\endce{\end{center}}
\def\begar{\begin{array}}
\def\endar{\end{array}}
\def\begeq{\begin{equation}}
\def\endeq{\end{equation}}
\def\begdi{\begin{displaymath}}
\def\enddi{\end{displaymath}}
\def\begdis{\begin{eqnarray*}}
\def\enddis{\end{eqnarray*}}
\def\begeqa{\begin{eqnarray}}
\def\endeqa{\end{eqnarray}}
\def\begdes{\begin{description}}
\def\enddes{\end{description}}
\def\begit{\begin{itemize}}
\def\endit{\end{itemize}}
\def\begen{\begin{enumerate}}
\def\enden{\end{enumerate}}
\def\beglar{\left[\begin{array}}
\def\endrar{\end{array}\right]}
\def\begle{\begin{lemma}}
\def\endle{\end{lemma}}
\def\begde{\begin{definition}}
\def\endde{\end{definition}}
\def\begth{\begin{theorem}}
\def\endth{\end{theorem}}
\def\begco{\begin{corollary}}
\def\endco{\end{corollary}}
\def\begprop{\begin{proposition}}
\def\endprop{\end{proposition}}	
\def\begex{\begin{example}}
\def\endex{\hfill\openbull \end{example}}
\def\begexer{\begin{exercise}}
\def\endexer{\end{exercise}}
\def\begalg{\begin{algo}}
\def\endalg{\end{algo}}
\def\begre{\noindent{\bf Remark}:\hspace*{0.05cm}}
\def\endre{\\}
\def\begres{\noindent{\bf Remarks}:\begin{enumerate}}
\def\endres{\end{enumerate} \par}
\def\begpr{\noindent{\em Proof:}$\;\;$}
\def\endpr{\hfill\bull}
\def\begtab{\begin{tabular}}
\def\endtab{\end{tabular}}
\def\rref#1{(\ref{#1})}
\newtheorem{lemma}{Lemma}[section]
\newtheorem{definition}{Definition}[section]
\newtheorem{theorem}{Theorem}[section]
\newtheorem{proposition}{Proposition}[section]
\newtheorem{corollary}{Corollary}[section]
\newtheorem{example}{Example}[section]
\newtheorem{algo}{Algorithm}[section]

\begin{document}
	\begin{frontmatter}
		
		\title{On Structural Non-commutativity in Affine Feedback of SISO Nonlinear Systems} 
		
		\thanks[footnoteinfo]{Supported by MASCOT project in UiT The Arctic Univeristy of Troms\o.}
		
		\author[First]{Venkatesh~G.~S.} 
		
		\address[First]{Department of Mathematics and Statistics, \\
			UiT The Arctic University of Troms\o, 9019 Troms\o, Norway.\\ (e-mail: subbarao.v.guggilam@uit.no)}
		
		\begin{abstract}                
		The affine feedback connection of SISO nonlinear systems modeled by Chen--Fliess series is shown to be a group action on the plant which is isomorphic to the semi-direct product of shuffle and additive group of non-commutative formal power series. The additive and multiplicative feedback loops in an affine feedback connection are thus proven to be structurally non-commutative. A flip in the order of these loops results in a net additive feedback loop.  
		\end{abstract}
		
		\begin{keyword}
		Nonlinear Systems, Affine Feedback, Structural Non-commutativity, Chen--Fliess Series, Formal Groups, Formal Lie Algebras.  
		\end{keyword}
		
	\end{frontmatter}
	
	\section{Introduction}\label{sec:Introduction}
\begin{figure}[h]
		\centering
		\includegraphics[scale=0.45]{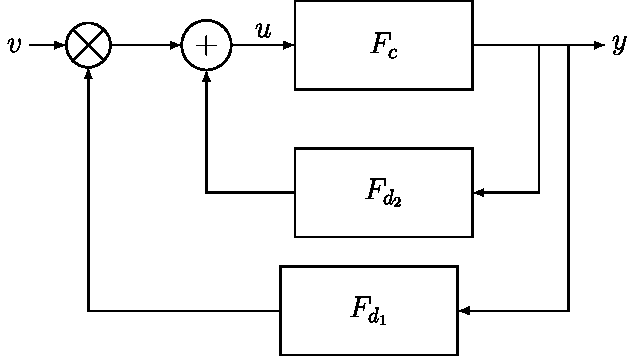}
		\caption{Affine Feedback Interconnection}
		\label{fig:affine-feedback}
	\end{figure}

	The goal of this paper is to {\em quantify} the structural non-commutativity between the additive and multiplicative feedback loops in an {\em affine feedback} connection of single-input single-output(SISO) nonlinear systems modeled by {\em Chen--Fliess} series. Chen--Fliess series is an weighted iterated integral series that models the input--output behavior of a (input-affine) nonlinear system. 	The prologue of using input-output maps for analysis of systems is that it is unique to the system (under analysis) and hence a coordinate free approach. Thus any properties proven via this paradigm is intrinsic  while the state space realizations are co-ordinate dependent. The argument is  akin to significance of Laplace-domain transfer functions/ time-domain input-output maps to state-space modeling in the domain of linear time-invariant systems. Affine feedback (as shown in Figure~\ref{fig:affine-feedback}) is a pivotal interconnection in the design of both static and dynamic feedback linearization of nonlinear systems~\citep{Isidori_95,Gray-Ebrahimi-Fard_SIAM17}.   
	
	The affine feedback interconnection of nonlinear systems modeled by Chen--Fliess was first dealt in the work of \cite{Gray-Ebrahimi-Fard_SIAM17}. The affine feedback group (transformation group in the affine feedback interconnection) was explicitly detailed there along with the Hopf algebra of coordinate functions on formal power series which provides a computational framework to compute the group product and the inverse. The work of \cite{Foissy_2016} has defined a {\em post-Lie} algebra purely from the Hopf algebra of coordinate functions corresponding to the affine feedback. Recently in the work of \cite{KEF_etal_2023}, the affine feedback group is proved to be the {\em Grossman--Larson} or {derived group} of a {\em post-group} (post group is the image of exponential map in the universal enveloping algebra of the post-Lie algebra). The document will not delve into post-groups and post-Lie algebras. For more details on post-groups and post-Lie algebras, refer~\cite{Alkaabi2023postgroup}. 
	
	The document shows firstly that affine feedback is a group action on the plant, where the group is same as the group in the definition of post-group (explained in Section~\ref{sec:affine-feedb-grp}). The commutator of two subgroups, corresponding to the multiplicative feedback loop and additive feedback loop in Figure~\ref{fig:affine-feedback} quantifies the structural non-commutativity between the loops viz., if one flips the order between the inner(additive) and outer(multiplicative) loops, the overall input--output map ($v \mapsto y$) will not remain the same. As such, proving the {\em existence} of structural non-commutativity is not a difficult question but the contribution of the paper stands out in {\em quantifying} it via the commutator of subgroups corresponding to additive and multiplicative feedback loops (in its Lie algebra at the infinitesimal level).    
	
	The organization of the document is as follows: Section~\ref{sec:formal-power-series} introduces notations and defines the shuffle algebra of formal power series. Section~\ref{sec:CF-series} gives a brief on Chen--Fliess series. Section~\ref{sec:affine-feedb-grp} defines the affine feedback transformation group and recollects some of the material from the works of \cite{Gray-Ebrahimi-Fard_SIAM17} and \cite{KEF_etal_2023}. Further, this section leads up to the main results of this document.  
	\section{Formal Power Series over $\re$} \label{sec:formal-power-series}
  The section recalls some preliminaries on non-commutative formal power series and bialgebras. For more details regarding the latter, we refer to~\citep{Sweedler}. 
	\smallskip
		
	\noindent Let $X = \{x_0,x_1\}$ denote an alphabet, i.e., a set of non-commuting indeterminates.
	\begen[(i)]
	\item $\allpoly$ is the $\re$-algebra of non-commutative polynomials with indeterminates in $X$. As a vector space, it is spanned by the set of words, denoted by $\allwords$, with letters from $X$. The unit is the empty word, which we write $\mathbf{1}$, and the set of non-empty words is denoted by $X^+:=X^*\backslash\{\mathbf{1}\}$.
	\smallskip
	\item $\allpoly$ is a cocommutative Hopf algebra, graded (by the word length, $|w| \ge 0$, $w \in X^*$) and connected, with catenation product and counit $\epsilon$ such that $\epsilon(w) = \delta_{|w|,0}$ for any $w \in X^*$. The coproduct is unshuffle $\Delta_{\shuffle}$ for which the elements $x_i \in X$ are primitive, i.e., $\Delta_{\shuffle}(x_i) = x_i \otimes \mathbf{1} + \mathbf{1} \otimes x_i$. A few examples of the unshuffle coproduct are given next. For letters $x_i,x_j \in X$:
	\begin{align*}
		\Delta_{\shuffle}x_ix_j 
		&= x_ix_j \otimes \mathbf{1} + x_i \otimes x_j + x_j \otimes x_i +  \\
		& \quad \quad \quad \mathbf{1} \otimes x_ix_j, \\
		\Delta_{\shuffle}x_i^2 
		&= x_i^2 \otimes \mathbf{1} + 2 x_i \otimes x_i + \mathbf{1} \otimes x_i^2. 
	\end{align*}
	\begre{
		Traditionally, the counit of $\allpoly$ is denoted by $\emptyset$ in the Chen--Fliess literature; for the rest of the article the presentation also follows this custom.
	}
	\smallskip
	\item The convolution algebra of linear maps from the counital coalgebra $\left(\allpoly, \Delta_{\shuffle}\right)$ to $\re$, is given by the space of formal power series denoted by $\allseries{}$. The dual basis is given by $\{\emptyset\} \cup X^{+}$, such that $\eta\left(\xi\right) = 1$ if $\eta = \xi$ in $X^+$, and zero else. An element $c \in \allseries{}$ is represented by
	\begin{align*}
		c = c(\mathbf{1})\emptyset + \sum\limits_{\eta \in X^+} c(\eta) \eta.
	\end{align*} 
	In the following, $\emptyset$ is not explicitly written unless needed. 
	\smallskip
	The convolution product on $\allseries{}$ is the shuffle product, which is defined for all $c,d \in \allseries{}$ and $p \in \allpoly$ by
	\begin{align*}
		\left(c \shuffle d\right)(p) 
		&= m_{\re} \circ \left(c \otimes d\right) \circ \Delta_{\shuffle} (p).
	\end{align*} 
	Here, $m_{\re}$ is the usual product in $\re$ and the unit element is the counit $\emptyset$. 
	\smallskip
	\item The maximal ideal of $\allseries{}$ is denoted by $\allproperseries{} : = \{c \in \allseries{} : c(\mbf{1}) = 0\}$. The unit group is denoted by $G_{\shuffle} : = \allseries{}\setminus\allproperseries{}$. The shuffle group $\left(M, \shuffle\right)$ is the normal subgroup of the unit group of $\allseries{}$ where 
	\begin{align*}
		M &= \{1 \emptyset + c': c' \in \allseries{}, c' \in \allproperseries{}\}.
	\end{align*}The group inverse in $M$ is defined as:
	\begin{align*}
		c^{\shuffle -1} = \sum_{k \geq 0} (\mbf{1}-c)^{\shuffle k}.
	\end{align*}
	\smallskip
	
	\item The additive group of formal power series, denoted by $G_{+}$ is defined on the set $\delta + \allseries{}$, where the group product $\diamond$ respectively inverse are defined by
	\begin{align*}
		(\delta + c) \diamond (\delta + d) &:=  \delta + (c + d),\\
		(\delta + c)^{\diamond -1} &:= \delta - c,
	\end{align*}
	for $c \in \allseries{}$. The unit element is $\delta$. The group $G_{+}$ can be interpreted as translation of additive group of formal power series $(\allseries{}, + , 0)$ by the ({\em fictitious}) element $\delta$ such that $F_{\delta}[u] = u$. The context of $\delta$ becomes more lucid with the document.
	\enden 
	Looking ahead, we define for $c \in \allseries{}$ the Chen--Fliess series $F_c[u](t)$ as a functional series such that the basis (dual words) of $\allseries{}$ are identified with iterated integrals of the function $u$, a control signal (detailed in Section~\ref{sec:CF-series}). In the context of a dynamical system, a Chen--Fliess series expresses its input-output behavior and provides a coordinate-independent framework to analyze the intrinsic properties of the system. 
	
	\section{Chen--Fliess Series}
	\label{sec:CF-series}
	
	The section acts as a review of Chen--Fliess series and the affine feedback interconnection. Let $u : [0,T] \longrightarrow \re$ be a $L_{1}$ measurable map over the compact interval $[0,T]$. Define the absolutely continuous function $U : [0,T] \longrightarrow \re$ such that $U(0) = 0$ as
	\begin{align*}
		U(t) = \int_{0}^t d\tau u(\tau).
	\end{align*}
	
	Define a monoid morphism $F$ from the monoid of (dual) words $\{\emptyset\} \cup X^{+}$ to the monoid of time-ordered iterated integrals of $u$ (where unit is constant map $1$) as
	\begin{align*}
		F_{x_1\eta}[u](t) &= \int_{0}^t dU_i(\tau)F_{\eta}[u](\tau),\\
		F_{x_0\eta}[u](t) &= \int_{0}^t d\tau F_{\eta}[u](\tau),
	\end{align*}
	where $\eta \in X^{+}\cup\{\emptyset\}$. The monoid morphism implies $F_{\emptyset}[u](t) = 1$.
	\begin{figure}[bt]
		\centering
		\includegraphics[scale= 0.65]{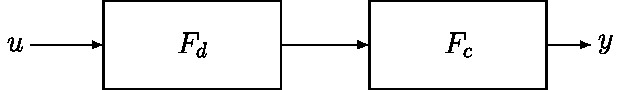}
		\caption{Composition of Chen--Fliess series}
		\label{fig:comp-CF-series}
	\end{figure}
	Given input $u: [0,T] \longrightarrow \re$ and a non-commutative formal power series $c$, the associated {\em Chen--Fliess} series, denoted by $F_{c}[u]$ is formal functional series defined as:
	\begin{align*}
		F_{c}[u](t) = c(\mbf{1}) + \sum_{\eta \in X^{+}} c(\eta)F_{\eta}[u](t).~ \quad \text{\citep{Fliess_81}}
	\end{align*}
	The map $F : c \mapsto F_c$ is an algebra homomorphism from $A_{\shuffle}$ to algebra of functions viz.~$(F_c.F_d)[u] = F_c[u]F_d[u] = F_{c\shuffle d}[u]$ for all $c,d \in \allseries{}$~\citep{Ree_58}. 
	
	The composition of Chen--Fliess series $F_c [u] \circ F_d [u]$, as illustrated in Figure~\ref{fig:comp-CF-series}, is a Chen--Fliess series $F_{c \circ d}[u]$. The formal power series $c \circ d$ is called the {composition product} of $c$ and $d$ and is defined as 
	\begin{align}\label{eqn:comp-prod}
		c \circ d = \sum_{\eta \in X^{*}} c\left(\eta\right) \eta \circ d,
	\end{align}
	where the composition product restricted to words is inductively defined as \citep{Ferfera_80}:
	\begin{align}\label{eqn:comp-prod-ind}
		\emptyset \circ d &= \emptyset, \nonumber \\
		x_0\eta \circ d  &= x_0\left(\eta \circ d\right), \\
		x_1\eta \circ d &= x_0\left(d \shuffle \left(\eta \circ d\right)\right) \nonumber,
	\end{align}
	for all $\eta \in X^{+} \cup \{\emptyset\}$. The following theorem states that composition product distributes over the shuffle product on the left. 
	\begth\citep{Ferfera_80}\label{th:comp-over-shuff} Let $c,d$ and $e \in \allseries{}$, then $\left(c\shuffle d\right) \circ e = \left(c \circ e\right) \shuffle \left(d \circ e\right)$.
	\endth
		
	Let $B_{p}(R,T) := \{u \in L_{p}[0,T]: \norm{u}_p < R\}$
	where $\norm{\cdot}_p$ is the usual $L_p$-norm. A Chen--Fliess $F_c$ is called {\em convergent} if is a well-defined mapping from $B_{1}(R,T)$ to $B_{\infty}(R,T)$ for $R, T > 0$ viz $F_{c}[u] \in L_{\infty}([0,T]; \re)$. For more details on convergence of Chen--Fliess series see \citep{Gray-Wang_02,Thitsa-Gray_SIAM12,GS_thesis}.	If convergence is guaranteed, a Chen--Fliess series $F_c$ is called a {\em Fliess operator}. A Fliess operator $F_c$ can be realized by a state-space realization if and only if its underlying formal power series $c$ has a finite Lie-rank~\citep{Fliess_81}. Thus, Fliess operators provide {\em more general} descriptors of an input-output relation of the underlying system. 

	\section{Affine Feedback Interconnection of Chen--Fliess Series}
	\label{sec:affine-feedb-grp}
	
	The goal of this section is to review and refurbish certain concepts of affine feedback in the context of Chen--Fliess series, as illustrated in Figure~\ref{fig:affine-feedback}, and defined in the work of~\cite{Gray-Ebrahimi-Fard_SIAM17}. Prior going to details of how the closed-loop affine feedback product is defined, in the context of Chen--Fliess series, it is necessary to review the affine feedback group.
	
	\subsection{Affine Feedback Group}\label{ssec:affine-feedb-grp}
	The affine feedback group is defined on a subset $G$ of two-dimensional vectors of formal power series, $\allseries{}^2$, expressed in the basis $\{\text{e}_1,\text{e}_2\}$. The set $G \subset \allseries{}^2$ is defined as
	\begin{align*}
		G := \{c_1\text{e}_1 + c_2\text{e}_2 : c_1(\mbf{1}) = 1\},
	\end{align*}
	whence $G \cong M \times \allseries{}$  as sets. The elements of $G$ are represented by bold-face font. Interchangeably, $\mathbf{c} = c_1\text{e}_1 + c_2\text{e}_2 \, \in G$ is represented as vector $[c_1 \;\; c_2]^t$, given the basis of $\allseries{}^2$ is fixed\footnote{Note that the basis vectors $\text{e}_1$ and $\text{e}_2$ are referred to as $\delta$ and $1$ in Section~\ref{sec:formal-power-series} as well as in the work of ~\cite{Gray-Ebrahimi-Fard_SIAM17}.}. The series $1e_1 + 0e_2$ in $G$ is denoted by $\mathbf{e}$.
	\medskip
	
	If $\mathbf{c} = [c_1 \;\; c_2]^t \in G$ then we define
	$$
	F_{\mathbf{c}}[u] = uF_{c_1}[u] + F_{c_2}[u].
	$$ 
	The component of $\text{e}_1$ is called {\em multiplicative component} (or channel) while the component of $\text{e}_2$ is termed as {\em additive component} (or channel). The composition of two formal power series, $\circ$ defined in Section~\ref{sec:CF-series} is extended, when its right argument is in the set $G$. 
	\begde\label{def:comp_G} Let $\mathbf{d} \in G$, with $\mathbf{d} = [d_1 \; \; d_2]^t$, and $c \in \allseries{}$, then the composition product of $\mathbf{d}$ and $c$. 
	\begin{align*}
		\mathbf{d} \bar{\circ} c = \begin{bmatrix}
			d_1 \circ c & d_2 \circ c
		\end{bmatrix}^t.
	\end{align*}
	\endde
	
\begin{figure}[tb]
		\includegraphics[scale=0.5]{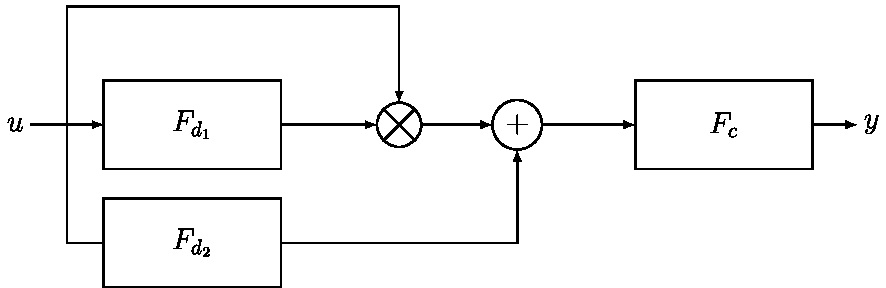}
		\caption{Composition of Chen--Fliess series:~$F_c$ and $uF_{d_1} + F_{d_2}$}
		\label{fig:mixed_comp}
	\end{figure}

	\begre~Given $\mathbf{d} \in G$ and $c \in \allseries{}$, it is to be noted that the composition of Chen--Fliess series $\left(F_{\mathbf{d}} \circ F_{c}\right)[u] \neq F_{\mathbf{d}\bar{\circ}c}[u]$. The $\bar{\circ}$ is merely extension of composition product, $\circ$, componentwise in $G$.
	\endre
	\smallskip
	The composition of Chen--Fliess series $F_{c} \circ F_{\mathbf{d}}$ where $c \in \allseries{}$ and $\mathbf{d} \in G$ (illustrated in Figure~\ref{fig:mixed_comp}), is given by the {\em mixed composition product} of $c$ and $\mathbf{d}$ giving a formal power series denoted by $c \modcomp \mathbf{d}$.
	
	\begde(Definition~3.1; \citep{Gray-Ebrahimi-Fard_SIAM17})\label{def:mix-pro-single}  Let $c \in \allseries{}$ and $\mathbf{d} = [d_1 \; d_2]^t \in G$, then the mixed composition product $c \modcomp \mathbf{d} \in \allseries{}$ is defined as
	\begin{align}\label{eqn:mixcomp-prod}
		c \modcomp \mathbf{d} = \sum_{\eta \in X^{*}} c(\eta) (\eta \modcomp \mathbf{d}),
	\end{align}
	where the mixed composition product (restricted to dual words) is defined inductively as
	\begin{align}\label{eqn:mix-comp-ind}
		\emptyset \modcomp \mathbf{d} &= \emptyset,\nonumber\\
		x_0\eta \modcomp \mathbf{d} &= x_0(\eta \modcomp \mathbf{d}), \\
		x_1\eta \modcomp \mathbf{d} &= x_1(d_1 \shuffle (\eta \modcomp \mathbf{d})) + x_0(d_2 \shuffle (\eta \modcomp \mathbf{d})),\nonumber
	\end{align}
	for all $\eta \in \{\emptyset\} \cup X^{+}$.
	\endde
	The mixed composition product is linear in its left argument (by definition) and distributes over the shuffle product, formally stated in the following theorem
	
	\begth\label{th:dist-over-sh}(Lemma~$3.3$,~\citep{Gray-Ebrahimi-Fard_SIAM17}) Let $c,d \in \allseries{}$ and $\mathbf{z} \in G$, then
	\begin{align*}
		\left(c \shuffle d\right) \modcomp \mathbf{z} = \left(c \modcomp \mathbf{z}\right) \shuffle \left(d \modcomp \mathbf{z}\right).
	\end{align*} 
	\endth 
	
	From~\rref{eqn:mix-comp-ind}, the following observations are evident and can be proved by induction on the length of words:
	\begin{enumerate}[(i)]\label{obv:modcomp}
		\item $x_0^k \modcomp \mathbf{d} = x_0^k$ for all $k \geq 0$ and $\mathbf{d} \in G$. 
		\item If $\mathbf{d} \neq \mathbf{e}$, $n \geq 1$ and $\eta = x_{i_1}x_{i_2}\cdots x_{i_n}$ such that $\exists j \in \{1,2,\ldots,n\}$ with $x_{i_j} = x_1$, then $\eta \modcomp \mathbf{d} \neq \eta$.    
		\item If for some $\mathbf{d} \in G$, $c \modcomp \mathbf{d} = c$ for all $c \in \allseries{}$, then $\mathbf{d} = \mathbf{e}$.
	\end{enumerate}
	
	Thus, if $\allseries{}_q$ := $\allseries{}/\sim$,  where the equivalence relation $\sim$ on $\allseries{}$ is defined as:
	\begin{align*}
		c\sim d \iff c-d \in \re[[x_0]] \quad \, \forall \, c,d \in \allseries{}.
	\end{align*}  
	
	\begle\label{lem:free-act-modcomp} If $c \in \allseries{}$, then its corresponding coset in $\allseries{}_q$ is denoted by $[c]$. Let $\mathbf{d} \in G$, then define $[c] \modcomp \mathbf{d} = [c \modcomp \mathbf{d}]$. The following statements are true: 
	
	\begen
	\item If $h \in [c]$ then $[h \modcomp \mathbf{d}] = [c \modcomp \mathbf{d}]$ (Mixed Composition product is well-defined on $\allseries{}_q$).
	\item If $[c] \neq [0]$, then $[c] \modcomp \mathbf{d} = [c]$ if and only if $\mathbf{d} = \mathbf{e}$. 
	\enden 
	\endle 
	Lemma~\ref{lem:free-act-modcomp} states that if the Chen--Fliess series $F_c$ in Figure~\ref{fig:mixed_comp} is affected by the input $u$ (having a least one word containing $x_1$ in the support of $c$), then $y = F_{c \modcomp \mathbf{d}} = F_{c}[u]$ if and only if $F_{\mathbf{d}}[u] = u$.  
	
	The following theorem asserts a sort of mixed associativity of composition product $\circ$ and mixed composition product $\modcomp$. 
	
	\begth\citep{KEF_etal_2023}\label{th:mix-assoc}Let $h,c \in \allseries{}$ and $\mathbf{d} \in G$, then $h \circ \left(c \modcomp \mathbf{d}\right) = \left(h \circ c\right) \modcomp \mathbf{d}$. 
	\endth
		
	The definition of mixed composition product, $\modcomp$, is extended as a set endomorphism of $G$ by applying the product componentwise. The extended product is denoted by $\triangleleft$ which is formally defined next.
	\begde \label{def:mix-pro}Let $\mathbf{c},\mathbf{d} \in G$ with $\mathbf{c} = [c_1\;\; c_2]^t$, then
	\begin{align*}
		\triangleleft: G \times G &\longrightarrow G \\
		(\mathbf{c},\mathbf{d}) &\mapsto \mathbf{c} \triangleleft\mathbf{d} := [c_1\modcomp \mathbf{d} \;\; c_2\modcomp \mathbf{d}]^t.
	\end{align*}
	\endde
	
	For Definition~\ref{def:mix-pro} to be valid, one needs to prove that $(c_1 \modcomp \mathbf{d})(\mbf{1}) = 1$ which is evident from \rref{eqn:mix-comp-ind}. Both $\triangleleft$ and $\modcomp$ shall be addressed as mixed composition product in the document with no loss of clarity as the domain of both maps are not identical.
	
	\begre~Given $\mathbf{c},\mathbf{d} \in G$, it is to be noted that the composition of Chen--Fliess series $\left(F_{\mathbf{c}} \circ F_{\mathbf{d}}\right)[u] \neq F_{\mathbf{c}\triangleleft\mathbf{d}}[u]$. The $\triangleleft$ is merely extension of the mixed composition product, $\modcomp$, componentwise in $G$. The product on $G$ encoding the composition $\left(F_{\mathbf{c}} \circ F_{\mathbf{d}}\right)[u]$ is in definition~\ref{def:circ_a_prod}.
	\endre
	
	The following theorem asserts the interplay between the extended composition product $\bar{\circ}$ and modified composition products $\modcomp$ and $\triangleleft$.
	\begth\label{th:bcir-mdc-lhk} Let $\mathbf{d}, \mathbf{y} \in G$, with $\mathbf{d} = \left[d_1 \;\; d_2\right]^t$, and $c \in \allseries{}$, then
	\begin{align*}
		\mathbf{d} \bar{\circ} \left(c \modcomp \mathbf{y}\right) = \left(\mathbf{d}\bar{\circ} c\right)\triangleleft \mathbf{y}.
	\end{align*} 
	\endth
	
	\begpr
	\begin{align*}
		\mathbf{d} \bar{\circ} \left(c \modcomp \mathbf{y}\right) &= \begin{bmatrix}
			d_1 \circ \left(c \modcomp \mathbf{y}\right) \\
			d_2 \circ \left(c \modcomp \mathbf{y}\right)
		\end{bmatrix}.
	\end{align*} 
	Using Theorem~\ref{th:mix-assoc} on each component, 
	\begin{align*}
		\mathbf{d} \bar{\circ} \left(c \modcomp \mathbf{y}\right) &= \begin{bmatrix}
			\left(d_1 \circ c \right) \modcomp \mathbf{y} \\
			\left(d_2 \circ c \right) \modcomp \mathbf{y} 
		\end{bmatrix} = \left(\mathbf{d} \bar{\circ} c\right) \triangleleft \mathbf{y}.
	\end{align*}
	\endpr
	\smallskip
	
	The set $G$ can be endowed with two different group structures that are vital in affine feedback interconnection, the first of which is explained next.
	
	\begth\citep{KEF_etal_2023}\label{th:odot-grp} Let $\mathbf{c},\mathbf{d} \in G$ with $\mathbf{c} = [c_1 \;\; c_2]^t$ and $\mathbf{d} = [d_1 \;\; d_2]^t$. Define a binary product $\odot$ on $G$ as\footnote{In \cite{KEF_etal_2023}, the notation is $\cdot$ instead of $\odot$.}
	\begin{align} \label{eqn:grp-prod-odot}
		\mathbf{c} \odot \mathbf{d} = 	\begin{bmatrix}
			c_1 \shuffle d_1 & c_2 + (d_2 \shuffle c_1)
		\end{bmatrix}^t,
	\end{align}
	then $(G, \odot)$ is a group with identity element being $\mathbf{e}:=[1 \;\; 0]^t$. The group inverse for $\mathbf{c} \in G$ is defined as
	\begin{align}\label{eqn:grp-inv-odot}
		\mathbf{c}^{\tiny{\odot-1}} = 	\begin{bmatrix}
			c_1^{\shuffle -1} & -c_{1}^{\shuffle -1} \shuffle c_2
		\end{bmatrix}^t.
	\end{align}
	\endth
	
	
	The group $(G, \odot)$ is isomorphic to the semi-direct product of $M$ and the additive group of formal power series.
	\begin{align}
		(G, \odot) \cong (M, \shuffle) \ltimes (\allseries{},+),
	\end{align}
	with the split exact sequence of groups given by
	\begeq\label{eq:exact-sequence}
	\begin{tikzcd}
		0 \arrow{r} & \allseries{} \arrow{r}{i} & G \arrow{r}{\pi_1} & M \arrow{r} & 1 \, ,
	\end{tikzcd}
	\endeq
	where the monomorphism $i(c) = 1e_1 + ce_2$ for all $ c \in \allseries{}$ and  $\pi_1$ is the canonical projection onto the first component. Observe that the additive group $G_{+} = \left(\delta + \allseries{}, \diamond\right)$ defined in Section~\ref{sec:formal-power-series} is the embedding of $(\allseries{}, +)$ in $G$. 
	
	\begth\citep{KEF_etal_2023} \label{thm:dist-over-odot}Let $\mathbf{x}, \mathbf{y}, \mathbf{z}$ in $G$, then 
	\begin{enumerate}[(i)]
		\item $\left(\mathbf{x}\odot\mathbf{y}\right) \triangleleft \mathbf{d} = \left(\mathbf{x} \triangleleft \mathbf{d}\right) \odot \left(\mathbf{y} \triangleleft \mathbf{d}\right)$.
		\item $\mathbf{e} \triangleleft \mathbf{d} = \mathbf{e}$,
	\end{enumerate}
	where $\mathbf{e}$ is the identity element of $\left(G, \odot\right)$.
	\endth
	
	%
	
	Theorem~\ref{thm:dist-over-odot} infers that $\triangleleft$ respects the group structure $\left(G,\odot\right)$ in its left argument. Identical to $\triangleleft$, the extended composition product, $\bar{\circ}$, respects the group structure in its left argument and is stated in the following theorem. 
	\begth\label{th:bcirc-over-odot} Let $\mathbf{x},\mathbf{y}$ be elements of the group $\left(G,\odot\right)$ and $c \in \allseries{}$, then
	\begin{enumerate}[(i)]
		\item $\left(\mathbf{x}\odot\mathbf{y}\right)\bar{\circ} c = \left(\mathbf{x} \bar{\circ} 0  c\right) \odot \left(\mathbf{y} \bar{\circ} c\right)$,
		\item $\mathbf{e} \bar{\circ} c = \mathbf{e}$,
	\end{enumerate}
	where $\mathbf{e}$ is the identity element of the group $\left(G,\odot\right)$.
	\endth
	
	\begpr (ii) is a straightforward verification from definition~\ref{def:comp_G}.
	(i) Suppose $\mathbf{x} = x_1e_1 + x_2e_2$ and $\mathbf{y} = y_1e_1 + y_2e_2$, then using \rref{eqn:grp-prod-odot} and definition~\ref{def:comp_G},
	\begin{align*}
		\left(\mathbf{x} \odot \mathbf{y}\right) \bar{\circ} c &= \begin{bmatrix}
			x_1 \shuffle y_1 \\ x_2 + \left(y_2 \shuffle x_1\right)
		\end{bmatrix} \bar{\circ} c \\
		&= \begin{bmatrix}
			\left(x_1 \shuffle y_1\right) \circ c \\ 
			\left(x_2 + \left(y_2 \shuffle x_1\right)\right) \circ c
		\end{bmatrix}.
	\end{align*}
	Using Theorem~\ref{th:comp-over-shuff} and \rref{eqn:grp-prod-odot},
	\begin{align*}
		\left(\mathbf{x} \odot \mathbf{y}\right) \bar{\circ} c &= \begin{bmatrix}
			\left(x_1 \circ c\right) \shuffle \left(y_1 \circ c\right)\\
			\left(x_2 \circ c\right) + \left(\left(y_2 \circ c\right) + \left(x_1 \circ c\right)\right)
		\end{bmatrix} \\
		&= \begin{bmatrix}
			x_1 \circ c \\ x_2 \circ c
		\end{bmatrix} \odot \begin{bmatrix}
			y_1 \circ c \\ y_2 \circ c
		\end{bmatrix} = \left(\mathbf{x}\bar{\circ}c\right) \odot \left(\mathbf{y}\bar{\circ} c\right).
	\end{align*}
	\endpr
	
	\begre{Theorem~\ref{thm:dist-over-odot} and Theorem~\ref{th:bcirc-over-odot} infers that should a composition product on formal power series $\allseries{}$ distribute over shuffle product on the left, then its component extension to set $G$ respects the group structure $\left(G, \odot\right)$ in the left argument.}
	\endre
	\smallskip
	
	The set $G$ possess another group structure, called {\em affine feedback group} and the group product is denoted by $\star$ whose definition as a binary product on the set $G$ is stated next. The definition for the product $\star$ stems\footnote{In \citep{Gray-Ebrahimi-Fard_SIAM17}, the notation for affine feedback group product is $\circ$.} from encoding the composition of two Chen--Fliess series whose generating series are from the set $G$ viz. $(F_{\mathbf{c}} \circ F_{\mathbf{d}})[u] = F_{\mathbf{c}\star \mathbf{d}}[u]$.
	
	\begde\label{def:circ_a_prod} Let $\mathbf{c}, \mathbf{d}$ in $G$, with $\mathbf{c} \star \mathbf{d}$, assuming $\mathbf{c} = c_1e_1 + c_2e_2$ and $\mathbf{d} = d_1e_1 + d_2e_2$, then
	\begin{align*}
		\star : G \times G &\longrightarrow G \\
		\left(\mathbf{c},\mathbf{d}\right) &\mapsto \left(\mathbf{c} \triangleleft \mathbf{d}\right) \odot \mathbf{d}.
	\end{align*}
	Thus, the coordinates of the product are given as
	\begeq\label{eqn:grp-prod-affine-cord}
	\mathbf{c} \star \mathbf{d} = \begin{bmatrix}
		\left(c_1 \modcomp \mathbf{d}\right) \shuffle d_1 & \left(\left(c_{1}\modcomp \mathbf{d}\right) \shuffle d_2\right) + \left(c_2\modcomp \mathbf{d}\right)  
	\end{bmatrix}^t.
	\endeq
	\endde
	
	\begth\label{thm:affine-feedb-grp}(Theorem~3.8;~\citep{Gray-Ebrahimi-Fard_SIAM17}) The tuple $\left(G,\star, \mathbf{e}\right)$ is a group where $\mathbf{e} := [1\;\; 0]^t$ is the identity element.
	\endth
	The group inverse of $\mathbf{c} \in (G, \star)$ is denoted by $\mathbf{c}^{\tiny{\star -1}}$. Given $\mathbf{c}$, $\mathbf{c}^{\tiny{\star -1}}$ can be computed by the antipode of Hopf algebra of coordinate maps on $\allseries{}$~\citep{Gray-Ebrahimi-Fard_SIAM17,Foissy_2016}. 
	
	\begth\label{thm:mix-comp-single-grp-act}(Theorem~$3.9$,~\citep{Gray-Ebrahimi-Fard_SIAM17}) Let $c \in \allseries{}$ and $\mathbf{x}, \mathbf{y} \in \left(G, \star\right)$, then
	\begin{enumerate}[(i)]
		\item $\left(c \modcomp \mathbf{x}\right)\modcomp \mathbf{y} = c \modcomp \left(\mathbf{x} \star \mathbf{y}\right)$,
		\item $c \modcomp \mathbf{e} = c$. 
	\end{enumerate}
	\endth
	Theorem~\ref{thm:mix-comp-single-grp-act} asserts that $\modcomp$ is a right action on the set of formal power series $\allseries{}$ by the affine feedback group $(G, \star)$. An important but evidently straightforward observation from  Theorem~\ref{thm:mix-comp-single-grp-act} and Lemma~\ref{lem:free-act-modcomp} is the following theorem.
	
	\begth\label{thm:freeness-mix-comp-sing} The following statements are true: 
	\begin{enumerate}
		\item $\modcomp$ is an effective but not a free right group action of affine feedback group $\left(G, \star\right)$ on the set $\allseries{}$.
		\item $\modcomp$ is a free right group action of feedback group $\left(G, \star\right)$ on the quotient set $\allseries{}_q\setminus{[0]}$.
	\end{enumerate}
	\endth
	
	Following Theorem~\ref{thm:mix-comp-single-grp-act}, Theorem~\ref{thm:dist-over-odot} and definition~\ref{def:mix-pro}, the following theorem is evident stating that mixed composition product , $\triangleleft$, is a right group action of affine feedback group $\left(G, \star\right)$ on the group $\left(G, \odot\right)$.
	
	\begth\citep{KEF_etal_2023}\label{th:lhook-right-act} Let $\mathbf{x}$ in the set $G$ and $\mathbf{y},\mathbf{z}$ in the group $\left(G, \star\right)$, then
	\begin{align*}
		\left(\mathbf{x}\triangleleft\mathbf{y}\right)\triangleleft \mathbf{z} &= \mathbf{x}\triangleleft\left(\mathbf{y}\star \mathbf{z}\right) \\
		\mathbf{x}\triangleleft \mathbf{e} &= \mathbf{x},
	\end{align*} 
	where $\mathbf{e}$ is the identity element of the group $\left(G,\star\right)$.
	\endth
	
	Theorem~\ref{thm:freeness-mix-comp-sing} asserts that the group action $\triangleleft$ is effective but not free on $(G, \odot)$. A simple counterexample to the freeness of the group action is the following: Let $\mathbf{c} = [1+x_0 \; \; x_0]^t \, \in \left(G, \odot\right)$. Then, $\mathbf{c} \triangleleft \mathbf{d} = \mathbf{c}$ for all $\mathbf{d} \in \left(G,\star\right)$. The stabilizer subgroup of $c$ being the whole of $\left(G, \star\right)$ implies that the right action $\triangleleft$ is not free.
	
	\begre{An interesting question would be to find a/the (or prove the non-existence of) non-trivial normal subgroup $H$ of $\left(G,\odot\right)$ such that $\triangleleft$ acts freely on the quotient group $\left(G/H, \odot\right)$ by the group $\left(G, \star\right)$. However, the question is deferred for the future work.}  
	\endre	
	
	To summarize this section, the set $G$ has the following interacting {\em double} group structure : $\left(G, \odot, \star, \triangleleft\right)$, where $\left(G, \odot\right)$ and $\left(G, \star\right)$ are groups for which $\mathbf{e} = 1e_1 + 0e_2$ is the identity element for both the groups while $\triangleleft$ is a right action on the group $\left(G, \odot\right)$ by the group $\left(G,\star\right)$.
	
	\begre{$\left(G,\odot,\triangleleft\right)$ actually forms a post-group with the affine feedback group $(G, \star)$ as its Grossman--Larson group. For more on this, refer the work of \cite{KEF_etal_2023}.}
	\endre	
	
	\subsection{Affine Feedback Interconnection} 
\begin{figure}[tb]
		\centering
		\includegraphics[scale=0.5]{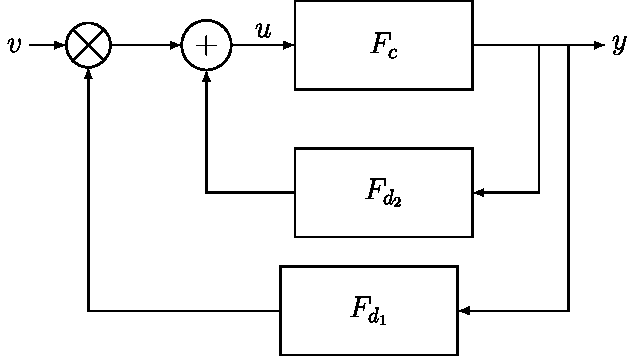}
		\caption{Affine Feedback Interconnection of $F_{c}$ with $F_{\mathbf{d}}$}
		\label{fig:affine-feedback2}
	\end{figure}

The {\em affine feedback product} describing the closed-loop of the affine feedback interconnection (as illustrated in Figure~\ref{fig:affine-feedback2}) is explained with the {\em newer} group structures endowed with interaction (called post-group) $(G, \odot, \triangleleft, \star)$. The affine feedback is shown to be a group action on plant and the {\em structural} non-commutativity between the additive and multiplicative feedback loops is proven and {\em quantified}. 
	
	\begre{The affine feedback product was first worked out in~\citep{Gray-Ebrahimi-Fard_SIAM17}. However the group action describing the closed-loop of the affine feedback and the structural non-commutativity of feedback loops are not evident from the cited work.}
	\endre 
	
	\begth \label{thm:affine-feedb}Let $c \in \allseries{}$ and $\mathbf{d} = d_1e_1 + d_2e_2 \in \left(G, \odot, \triangleleft, \star\right)$. The affine feedback interconnection of the Chen--Fliess series $F_c$ with the Chen--Fliess series $F_{\mathbf{d}}$, as in Figure~\ref{fig:affine-feedback2}, is represented by a Chen--Fliess series $F_{c @\mathbf{d}}$, where the affine feedback product $c@\mathbf{d} \in \allseries{}$ is defined as
	\begin{align}\label{eqn:aff-feedb}
		c @ \mathbf{d} = c \modcomp \left(\left(\mathbf{d}^{\odot -1}\right) \bar{\circ} c\right)^{\tiny{\star -1}}.
	\end{align}
	\endth
	
	\begpr Taking the expression on the RHS and using Theorem~\ref{th:bcirc-over-odot} and \rref{eqn:grp-inv-odot},
	\begin{align*}
		\left(\mathbf{d}^{\tiny{\odot -1}}\right)\bar{\circ} c &= \left(\mathbf{d} \bar{\circ} c\right)^{\tiny{\odot -1}} \\
		&= \begin{bmatrix}
			\left(d_1 \circ c \right)^{\shuffle -1} \\
			-\left(d_1 \circ c\right)^{\shuffle -1} \shuffle \left(d_2 \circ c\right)
		\end{bmatrix}.
	\end{align*}
	Thus, 
	\begin{align*}
		c @ \mathbf{d} &= c \modcomp \left(\left(\mathbf{d}^{\odot -1}\right) \bar{\circ} c\right)^{\tiny{\star -1}}\\
		&= c \modcomp \begin{bmatrix}
			\left(d_1 \circ c \right)^{\shuffle -1} \\
			-\left(d_1 \circ c\right)^{\shuffle -1} \shuffle \left(d_2 \circ c\right)
		\end{bmatrix}^{\tiny{\star -1}}.
	\end{align*}  
	The expression on the RHS is exactly the affine feedback product defined in $(3.7)$ of ~\citep{Gray-Ebrahimi-Fard_SIAM17}.
	
	\endpr
	\medskip
	
	The notion that feedback can described mathematically as a transformation group acting on the plant is well established in control theory~\citep{Brockett_78}. The following theorem asserts the group action in the context of affine feedback as illustrated in Figure~\ref{fig:affine-feedback2}.
	
	\begth\label{th:feedb-grp-action} Let $c \in \allseries{}$ and $\mathbf{x}, \mathbf{y} \in \left(G, \odot, \triangleleft, \star \right)$, then
	\begeqa\label{eqn:feedb-grp-odot-action}
	\left(c @ \mathbf{x}\right)@ \mathbf{y} &=& c @ \left(\mathbf{x} \odot \mathbf{y}\right), \nonumber \\ 
	\left(c@\mathbf{e}\right) &=& c,
	\endeqa 
	where $e$ is the identity element of $\left(G,\odot, \triangleleft, \star\right)$. The affine feedback is thus a right action by the group $\left(G, \odot\right)$ on the set $\allseries{}$.
	\endth
	
	\begpr Using Theorem~\ref{thm:affine-feedb},
	\begin{align*}
		\left(c @ \mathbf{x}\right)@\mathbf{y} &= \left(\left(c @\mathbf{x}\right) \modcomp \left(\mathbf{y}\right)^{\tiny{\odot -1}} \bar{\circ} \left(c @ \mathbf{x}\right)\right)^{\tiny{\star-1}} \\
		&=\left[c \modcomp \left(\left(\mathbf{x}\right)^{\tiny{\odot -1}}\bar{\circ} c\right)^{\tiny{\star -1}}\right] \modcomp \\
		& \quad \quad  \left[\left(\mathbf{y}\right)^{\tiny{\odot -1}} \bar{\circ} \left(c \modcomp \left(\left(\mathbf{x}\right)^{\tiny{\odot -1}} \bar{\circ} c\right)^{\tiny{\star -1}}\right)\right]^{\tiny{\star -1}}.
	\end{align*}
	Using Theorem~\ref{th:bcir-mdc-lhk},
	\begin{align*}
		\left(c @ \mathbf{x}\right)@\mathbf{y} &= \left[c \modcomp \left(\left(\mathbf{x}\right)^{\tiny{\odot -1}}\bar{\circ} c\right)^{\tiny{\star -1}}\right] \modcomp \\
		& \left[\left(\left(\mathbf{y}\right)^{\tiny{\odot -1}} \bar{\circ} c\right) \triangleleft \left(\left(\mathbf{x}\right)^{\tiny{\odot}-1}\bar{\circ}c\right)^{\tiny{\star -1}}\right]^{\tiny{\star -1}}.
	\end{align*}
	Using Theorem~\ref{thm:mix-comp-single-grp-act} and using the fact $\mathbf{a}^{\tiny{\star -1}}\star\mathbf{b}^{\tiny{\star-1}} = \left(\mathbf{b}\star \mathbf{a}\right)^{\tiny{\star -1}}$ for all $\mathbf{a},\mathbf{b} \in \left(G, \star\right)$,
	\begin{align*}
		\left(c @ \mathbf{x}\right)@\mathbf{y} &= c \modcomp \Bigg[\left(\left(\mathbf{x}\right)^{\tiny{\odot -1}}\bar{\circ} c\right)^{\tiny{\star -1}} \star \\
		& \left[\left(\left(\mathbf{y}\right)^{\tiny{\odot -1}} \bar{\circ} c\right) \triangleleft \left(\left(\mathbf{x}\right)^{\tiny{\odot}-1}\bar{\circ}c\right)^{\tiny{\star -1}}\right]^{\tiny{\star -1}}\Bigg] 
	\end{align*}
	\begin{align*}
		&= c \modcomp \Bigg[\left(\left(\left(\mathbf{y}\right)^{\tiny{\odot -1}}\bar{\circ}c\right) \triangleleft \left(\left(\mathbf{x}\right)^{\tiny{\odot-1}}\bar{\circ}c\right)^{\tiny{\star-1}}\right) \star \\
		& \quad \quad \quad \quad  \left(\left(\mathbf{x}\right)^{\tiny{\odot -1}}\bar{\circ}c\right)\Bigg]^{\tiny{\star-1}}.
	\end{align*}
	Using definition~\ref{def:circ_a_prod},
	\begin{align*}
		\left(c @ \mathbf{x}\right)@\mathbf{y} &= c \modcomp \Bigg[\Big(\left(\left(\left(\mathbf{y}\right)^{\tiny{\odot -1}}\bar{\circ}c\right)\triangleleft 
		\left(\left(\mathbf{x}\right)^{\tiny{\odot -1}}\bar{\circ}c\right)^{\tiny{\star -1}}\right) \\ 
		& \triangleleft \left(\left(\mathbf{x}\right)^{\tiny{\odot-1}}\bar{\circ}c\right)\Big) \odot \left(\left(\mathbf{x}\right)^{\tiny{\odot -1}}\bar{\circ}c\right)\Bigg]^{\tiny{\star -1}}.
	\end{align*}
	Using Theorem~\ref{th:lhook-right-act},
	\begin{align*}
		\left(c @ \mathbf{x}\right)@\mathbf{y} &= c \modcomp \Bigg[\Big(\left(\left(\mathbf{y}\right)^{\tiny{\odot-1}}\bar{\circ}c\right) \triangleleft \\
		& \quad \left[\left(\left(\mathbf{x}\right)^{\tiny{\odot -1}}\bar{\circ}c\right)^{\tiny{\star-1}}\star\left(\left(\mathbf{x}\right)^{\tiny{\odot-1}}\bar{\circ}c\right)\right]\Big) \odot \\ 
		&\quad \quad \quad  \left(\left(\mathbf{x}\right)^{\tiny{\odot-1}}\bar{\circ}c\right)\Bigg]^{\tiny{\star-1}} \\
		&= c \modcomp \Bigg[\left(\left(\left(\mathbf{y}\right)^{\tiny{\odot-1}}\bar{\circ}c\right) \triangleleft \mathbf{e}\right) \odot \\
		& \quad \quad \quad  \left(\left(\mathbf{x}\right)^{\tiny{\odot-1}}\bar{\circ}c\right)\Bigg]^{\tiny{\star-1}},
	\end{align*}
	where $\mathbf{e}$ is the identity element of group $\left(G,\star\right)$. Using Theorem~\ref{th:lhook-right-act} and Theorem~\ref{th:bcirc-over-odot},
	\begin{align*}
		\left(c @ \mathbf{x}\right)@\mathbf{y} &= c \modcomp \left[\left(\left(\mathbf{y}\right)^{\tiny{\odot -1}}\bar{\circ}c\right) \odot \left(\left(\mathbf{x}\right)^{\tiny{\odot -1}}\bar{\circ}c\right)\right]^{\tiny{\star -1}} \\
		&= c \modcomp \left[\left(\mathbf{x}\odot \mathbf{y}\right)^{\tiny{\odot -1}}\bar{\circ} c\right]^{\tiny{\star -1}} = c@\left(\mathbf{x}\odot \mathbf{y}\right).
	\end{align*}
	The statement $c@\mathbf{e} = c$ can be verified straightforward from the definition.
	\endpr
	
	Thus, the affine feedback product is a right action by the group $\left(G,\odot\right)$ on the set of formal power series $\allseries{}$.
	\begin{figure}[tb]
		\centering
		\includegraphics[scale = 0.4]{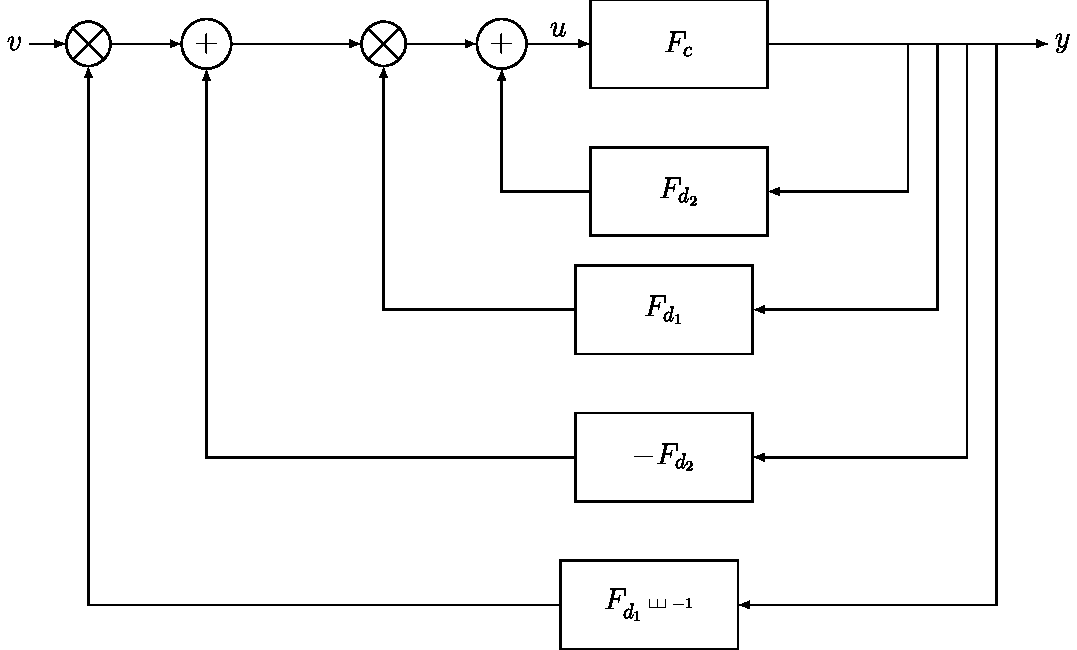}
		\caption{Non-commutativity of Additive and Multiplicative Feedback Loops.}
		\label{fig:4-loop}
	\end{figure}
	
	\begre{Additive feedback is an affine feedback whose multiplicative component is $1$, while multiplicative feedback is an affine feedback whose additive component is $0$. Thus, by \rref{eq:exact-sequence}, additive feedback product is action on the plant by the additive group of formal power series $\left(\allseries{}, +, 0\right)$ and respectively, the multiplicative feedback product is action on the plant by the shuffle group $(M, \shuffle, 1\emptyset)$. The assertions are in accordance with the literature~\citep{Gray-etal_SCL14,Venkat_MTNS_22,KEF_GS_2023}.}
	\endre
	
	The additive and multiplicative feedback loops in the affine feedback configuration do not commute in general. This structural non-commutativity is illustrated via the Figure~\ref{fig:4-loop}. Let $F_{\omega}[v] = y$, then the closed-loop generating series $\omega$ for the Figure~\ref{fig:4-loop} can be described as:
	\begin{align}\label{eq:nc-loop}
		\omega & = \left(\left(\left(c @ \begin{bmatrix}
			1 \\ d_2
		\end{bmatrix}\right)@\begin{bmatrix}
			d_1 \\ 0
		\end{bmatrix}\right)@\begin{bmatrix}
			1 \\ -d_2
		\end{bmatrix}\right)@\begin{bmatrix}
			d_{1}^{\shuffle -1} \\ 0
		\end{bmatrix}.
	\end{align}
	
	If multiplicative and additive feedback loops did commute, then the additive feedback loops would cancel other (since $d_2 + (-d_2) = 0$) and multiplicative loops would negate each other (since $d_1 \shuffle d_1^{\shuffle -1} = 1$). Hence, if the loops commuted, Figure~\ref{fig:4-loop} is equivalent to the following diagram.	
	\begin{figure}[h]
		\centering
		\includegraphics[scale=0.5]{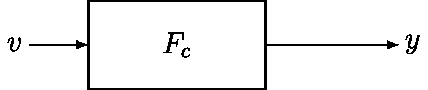}
	\end{figure}
	
\begin{figure}[h]
		\centering
		\includegraphics[scale=0.5]{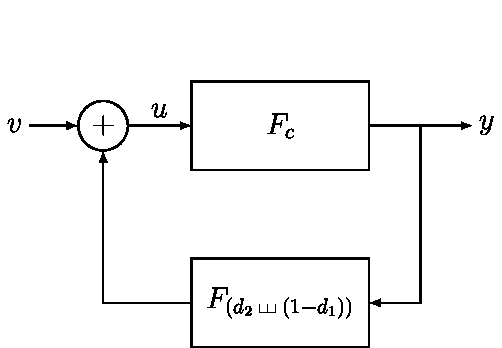}
		\caption{Residual Additive Feedback Loop of Figure~\ref{fig:4-loop}.}
		\label{fig:resd_loop}
	\end{figure}

	However, using Theorem~\ref{th:feedb-grp-action} on \rref{eq:nc-loop}:
	\begin{align*}
		\omega &= c @ \left(\begin{bmatrix}
			1 \\ d_2
		\end{bmatrix} \odot \begin{bmatrix}
			d_1 \\ 0
		\end{bmatrix} \odot \begin{bmatrix}
			1 \\ -d_2
		\end{bmatrix} \odot \begin{bmatrix}
			d_1^{\shuffle -1} \\ 0
		\end{bmatrix} \right).
	\end{align*}
	Using Theorem~\ref{th:odot-grp},
	\begin{align*}
		\omega &= c @ \left(\begin{bmatrix}
			d_1 \\ d_2
		\end{bmatrix} \odot \begin{bmatrix}
			d_1^{\shuffle -1} \\ -d_2
		\end{bmatrix}\right) \\
		&= c @ \begin{bmatrix}
			1 \\ d_2 \shuffle \left(1- d_1\right)
		\end{bmatrix}.
	\end{align*}
	Hence, $\omega = c@ \begin{bmatrix}
		1 & d_2 \shuffle \left(1-d_1\right)
	\end{bmatrix}^t$ infers that the feedback loops do not negate each other completely, and there is a {\em residual} additive feedback loop (since multiplicative component is $1$) amounting to $F_c$ being in additive feedback with the $F_{d_2 \shuffle (1-d_1)}$. Hence, the diagram~\ref{fig:4-loop} is equivalent to the Figure~\ref{fig:resd_loop}.
	
	Let $[\mathbf{x},\mathbf{y}]_{\odot}$ denote the group commutator of $\mathbf{x}, \mathbf{y}$ of $\left(G, \odot\right)$ viz.
	\begin{align*}
		[\mathbf{x}, \mathbf{y}]_{\odot} = \mathbf{x}\odot\mathbf{y}\odot\mathbf{x}^{\tiny{\odot-1}}\odot \mathbf{y}^{\tiny{\odot -1}}.
	\end{align*}
	
	If $A$ and $B$ are two subgroups of $\left(G, \odot\right)$, then
	\begin{align*}
		[A,B]_{\odot} := \{\mathbf{x}\odot\mathbf{y}\odot\mathbf{x}^{\tiny{\odot-1}}\odot\mathbf{y}^{\odot -1}: \, \forall \mathbf{x} \in A, \mathbf{y} \in B\}.
	\end{align*}
	
	The computation of $\omega$ in \rref{eq:nc-loop} asserts that $\left[M,G_{+}\right]_{\odot} \subseteq G_{+}$ and is formally proven in the following theorem.
	
	\begth \label{th:lower-series} The commutator subgroup $\widehat{\left[M,G_{+}\right]_{\odot}}$, defined as the group generated by the commutators $\left[M,G_{+}\right]_{\odot}$, is a normal subgroup of $G_{+}$. 
	\endth
	
	\begpr (i) Let $\mathbf{x} \in M$ and $\mathbf{y} \in G_{+}$, then
	\begin{align*}
		\left[\mathbf{x},\mathbf{y}\right]_{\odot} =  \mathbf{x}\odot\mathbf{y}\odot\mathbf{x}^{\tiny{\odot-1}}\odot\mathbf{y}^{\odot -1}.
	\end{align*}
	Since $G_{+}$ is a normal subgroup of $G$, this implies $\mathbf{x}\odot\mathbf{y}\odot\mathbf{x}^{\tiny{\odot}-1} \in G_{+}$.
	Thus,
	\begin{align*}
		\left[\mathbf{x},\mathbf{y}\right]_{\odot} =  \lbrace\mathbf{x}\odot\mathbf{y}\odot\mathbf{x}^{\tiny{\odot-1}}\rbrace\odot\mathbf{y}^{\odot -1} \in G_{+}
	\end{align*}
	Therefore, $\left[M,G_{+}\right]_{\odot} \subseteq G_{+} \Longrightarrow
	\widehat{\left[M,G_{+}\right]_{\odot}} \triangleleft G_{+}$ since a subgroup of a normal subgroup is normal.
	\endpr
	\smallskip	
	
	Let $\mathbf{c} \in M$ such that $\mathbf{c} = ce_1 + 0e_2$. Let $\mathbf{d} \in G_{+}$ such that 
	$\mathbf{d} = 1e_1 + de_2$. Then,
	\begin{align*}
		\left[\mathbf{b},\mathbf{d}\right]_{\odot} &= \begin{bmatrix}
			b \\ 0
		\end{bmatrix} \odot \begin{bmatrix}
			1 \\ d
		\end{bmatrix} \odot \begin{bmatrix}
			b^{\shuffle -1} \\ 0
		\end{bmatrix} \odot \begin{bmatrix}
			1 \\ -d
		\end{bmatrix} \\
		&= \begin{bmatrix}
			1 \\ d \shuffle \left(1-b\right)
		\end{bmatrix}.
	\end{align*}
	Hence, as a definition
	\begin{align}\label{eq:commutator}
		\left[M,G_{+}\right]_{\odot} := \Big\{1e_1 + \left(d \shuffle \left(1-b\right)\right):& d \in \allseries{}, \nonumber\\
		\quad \quad \quad b \in M\Big\}.
	\end{align}
	
	$\left(G, \odot\right)$ isomorphic to the character group of the graded connected combinatorial Hopf algebra of coordinate functions as discussed in~\citep{KEF_etal_2023}. $(G, \odot) \cong M \ltimes (\allseries{},+)$ is thus a formal Lie group. 
	The formal Lie algebra of $G$ is denoted by $\left(\mathfrak{g}, [\cdot,\cdot]\right)$. The Lie algebras of formal Lie subgroups $M$ and $G_{+}$ are denoted $\mathfrak{g}_{\shuffle}$ and $\mathfrak{g}_{+}$ respectively.
	
	Observe that $g_{\shuffle}$ is isomorphic to the vector space of $\allproperseries{}$. If for some $c \in \allseries{}$, $\mathbf{e} + tce_{1} \in M$ for $t\in \re$ implies 
	\begin{align*}
		t^{-1}\left(\left(\left(\mathbf{e} + tce_1\right) - \mathbf{e}\right)\left(\mathbf{1}\right)\right) = 0
		\Leftrightarrow c\left(\mathbf{1}\right) &= 0.
	\end{align*}
	Thus $c$ is a proper series.
Similarly, $\mathfrak{g}_{+}\cong \allseries{}$ as vector spaces. If for some $c \in \allseries{}$, $1e_1  + tce_2 \in (G_{+}, \diamond)$ for $t \in \re$ implies
	\begin{align*}
		t^{-1}\left((\mathbf{e}+ tce_2)- \mathbf{e}\right) = 0e_1 + ce_2.
	\end{align*}
	Therefore, $\mathfrak{g}_{+} \cong \allseries{}$ as vector spaces. 
	
	The commutator in \rref{eq:commutator} is  a map on the formal Lie groups given by:
	\begin{align*}
		\left[\cdot,\cdot\right]_{\odot} : M \times G_{+} &\longrightarrow  G_{+} \\
		\left(ce_1+ 0e_2, 1e_1+de_2\right) &\mapsto 1e_1 + \left(d \shuffle\left(1-c\right)\right)e_2,
	\end{align*}
	where $c \in M$ and $d \in \allseries{}$. Linearizing the group commutator map at $\mathbf{e}$, the identity of $\left(G, \odot\right)$, gives the Lie bracket between the Lie subalgebras $\mathfrak{g}_{\shuffle}$ and $\mathfrak{g}_{+}$. Let $d_1e_1 \in \mathfrak{g}_{\shuffle}$ and $d_2e_2 \in \mathfrak{g}_{+}$, where $d_1 \in \allproperseries{}, d_2 \in \allseries{}$. Then using \rref{eq:commutator},
	\begin{align*}
		\left[d_1e_1,d_2e_2\right] &= \left([(1+d_1)e_1,d_2e_2]_{\odot}-[1e_1,0e_2]_{\odot}\right) \\
		&= \begin{bmatrix}
			1  \\ d_2 \shuffle d_1  
		\end{bmatrix}  - \begin{bmatrix}
			1 \\ 0
		\end{bmatrix} = \begin{bmatrix}
			0 \\ d_1 \shuffle d_2
		\end{bmatrix}.
	\end{align*}
	Therefore the formal Lie group $\left(G,\odot\right)$ has a formal Lie algebra $\mathfrak{g}\cong \allproperseries{} \times \allseries{}$ (as vector spaces), whose Lie bracket is defined as the following
	\begin{align} \label{eq:Lie-brakcet}
		[d_1e_1,d_2e_2] = -[d_2e_2,d_1e_1] &= (d_1\shuffle d_2)e_2 \\
		[d_1e_1,f_1e_1] = [d_2e_2,f_2e_2] &= 0 \nonumber,
	\end{align} 
	where $d_1,f_1 \in \allproperseries{}$ and $d_2,f_2 \in \allseries{}$. Note that $\mathfrak{g}_{\shuffle}$ and $\mathfrak{g}_{+}$ are Abelian Lie algebras (as the corresponding group products are commutative). Henceforth in \rref{eq:Lie-brakcet}, 
	\begin{align*}
		[d_ie_i,f_ie_i] = 0 \quad \text{for} \; i = 1,2.
	\end{align*}
	
	\begre{
		The Lie algebra $\left(\mathfrak{g},[\cdot,\cdot]\right)$ in \rref{eq:Lie-brakcet} is the Lie algebra which is a part of definition of the post-Lie algebra that is obtained by linearizing the affine feedback post group $\left(G, \odot, \triangleleft\right)$, in \cite{KEF_etal_2023}. The post-Lie algebra was however first observed from the Hopf algebra of coordinate functions corresponding to affine feedback~\citep{Foissy_2016}.}
		\endre
		
	
	\begin{ack}
		The author acknowledges the support and various useful comments from  Prof.~Kurusch Ebrahimi-Fard of Department of Mathematical Sciences, Norwegian University of Science and Technology~(NTNU), Trondheim, Norway and Center for Advanced Studies~(CAS), Oslo, Norway and Prof. W. Steven Gray of Department of Electrical and Computer Engineering, Old Dominion University~(ODU), Norfolk, Virginia, USA. 	
	\end{ack}

\end{document}